\documentclass[12pt]{amsart}
\usepackage{amsmath,latexsym, amssymb}
\usepackage[all]{xy}
\usepackage[usenames]{color}
\usepackage[normalem]{ulem}
\def\to{\longrightarrow}
\def\ol{\overline}
\def\inv{^{-1}}
\def\id{\operatorname{id}}
\def\Tor{\operatorname{Tor}_1^B}
\def\coker{\operatorname{coker}}
\def\M{\mathcal{M}}
\def\Cbar{{\overline{C}}}

\def\CC{{\mathcal C}}
\def\FF{\mathcal{F}}
\def\PP{\mathcal{P}}
\def\EE{\mathcal{E}}
\def\GG{\mathcal{G}}
\def\RR{{\mathcal R}}

\def\e{{\varepsilon}}

\newcommand\Coa[1]{{\mbox{{\bf Coalg}-}#1}}
\newcommand\C[1]{\M_#1}
\newcommand\CMod[1]{\M od_#1}

\newtheorem{thm}{Theorem}[section]
\newtheorem{cor}[thm]{Corollary}
\newtheorem{prop}[thm]{Proposition}
\newtheorem{lem}[thm]{Lemma}

\theoremstyle{definition}
\newtheorem{defn}[thm]{Definition}
\newtheorem{example}[thm]{Example}

\title{On Radicals of Module Coalgebras}
\author{Yuqun Chen}
\thanks{The first author is supported by NNSF of China (Grant No. 10471045)
and NSF of Guangdong Province (Grant No. 021073 and 06025062).}
\address{School of Mathematical Sciences, South China Normal University, Guangzhou 510631, PRC.}
\email{yqchen@scnu.edu.cn}
\author{Siu-Hung Ng}
 \address{Department of Mathematics, Iowa State University, Ames, IA 50011, USA.}
 \email{rng@iastate.edu}
 \thanks{The second author is supported by NSA grant no. H98230-05-1-0020.}
 \author{Kar-Ping Shum}
 \address{Faculty of Science, The Chinese University of Hong Kong,Hong Kong,China (SAR).}
 \email{kpshum@math.cuhk.edu.hk}
 \thanks{The third author is partially supported by a RGC (HK) grant no. 2060297 (2005/07)}
\begin{document}
\begin{abstract}
We introduce the notion of idempotent radical class  of
module coalgebras over a bialgebra $B$. We prove that if $\RR$ is an
idempotent radical class of $B$-module coalgebras, then every
$B$-module coalgebra contains a unique maximal $B$-submodule
coalgebra in $\RR$. Moreover, a $B$-module coalgebra $C$ is a member
of $\RR$ if, and only if, $DB$ is in $\RR$ for every simple
subcoalgebra $D$ of $C$. The collection of $B$-cocleft coalgebras,
and the collection of $H$-projective module coalgebras over a Hopf
algebra $H$ are idempotent radical classes. As applications, we use
these idempotent radical classes to give another proofs for a
projectivity theorem and a normal basis theorem of Schneider without
assuming bijective antipode.
\end{abstract}
\maketitle
\section{Introduction}
In \cite{Ng98}, it has been proved that for any module coalgebra $C$
over a finite-dimensional Hopf algebra $H$, there exists a unique
maximal $H$-submodule coalgebra $\PP(C)$ of $C$ such that $\PP(C)$ is a projective $H$-module.
More
importantly, the ($H$-)projectivity of $C$ can be completely
determined by the projectivity of the $H$-submodules generated by
the simple subcoalgebras of $C$. Projectivity of modules over any
ring, in general, is not preserved under direct limit. However,
projectivity of $H$-module coalgebras behaves quite differently. Indeed,
the collection of projective $H$-module coalgebras carries certain
properties which allow the existence of unique
maximal projective $H$-submodule subcoalgebra of any given $H$-module coalgebra. This example is a
prototype of idempotent radical classes introduced in the sequel.\\

In this paper, we will give the definitions of {\em radical class}
and {\em idempotent radical class} of the category $\Coa{B}$  of
right module coalgebras over a bialgebra $B$. They are indeed
generalizations of coradical of coalgebras. We show that if $\RR$ is
an idempotent radical class of $\Coa{B}$, then for any $C \in
\Coa{B}$, there exists a unique maximal subcoalgebra $\RR(C)$ of $C$
such that $\RR(C)$ is a member of $\RR$. Moreover, $C \in \RR$ if,
and only if, $DB \in \RR$ for all simple subcoalgebras
$D$ of $C$.\\

To demonstrate the abundance  of idempotent radical
classes, we shall prove that the following well-known classes of right
$B$-module coalgebras are idempotent radical classes:
\begin{enumerate}
\item[(i)]  the collection $\PP$ of all projective $H$-module coalgebras over a Hopf algebra $H$;
\item[(ii)] the collection $\EE$ of all $B$-module coalgebras $C$ for which the functor
            $$
            \ol ?: \C{B}^C \to \M^{\ol C},\quad M \mapsto \ol M = M/MB^+
            $$
            defines an equivalence;
\item[(iii)] the collection $\CC$ of all cocleft $B$-module coalgebras.
\end{enumerate}
As applications of these idempotent radical classes, we shall give another proofs for a projectivity theorem and
a normal basis theorem of Schneider \cite{Schn90}. Obviously, Radford's freeness
theorem \cite{Radf77} for pointed Hopf algebras is an immediate consequence of Schneider's results.\\

Throughout this paper, we will assume all the algebras and
coalgebras are over the same ground field $k$ unless stated
otherwise. The tensor product $\otimes_k$ over the base field $k$ will
simply be denoted by $\otimes$, and  we always use $B$ to denote a bialgebra over $k$.

\section{Radical classes of Module Coalgebras}
Let $B$ be a bialgebra over a field $k$.
A coalgebra $C$ over $k$ is
said to be a (right) {\em $B$-module coalgebra}, or a {\em module
coalgebra over $B$} if $C$ is a right $B$-module and the coalgebra structure
maps $\Delta:C \to C \otimes C$ and $\e: C \to k$ are $B$-module maps, where
$C \otimes C$ is the right $B$-module with diagonal $B$-action and $k$ is considered as a
trivial $B$-module. The category of all (right) $B$-module coalgebras
will be denoted by $\Coa{B}$.\\

For any $B$-module coalgebra $C$, we
 will simply call a subcoalgebra of $C$, which is invariant under the
$B$-action, a {\em $B$-submodule coalgebra} of $C$.\\

If $X, Y$ are subspaces of a $B$-module coalgebra $C$, recall \cite{Sw69} that
the ``wedge" $X \wedge Y$ is defined as
$$
X \wedge Y= \Delta\inv(X \otimes C + C \otimes Y)\,.
$$
Following \cite{Mont93bk} and \cite{Sw69}, we define $\bigwedge^1 X
= X$ and $\bigwedge^{n+1}X =X \wedge (\bigwedge^n X)$ for $n \ge 1$.
If $X$ contains the coradical $C_0$ of $C$, then we have
$$
\sum_{n\ge 1} \bigwedge^n X = C\,.
$$
In particular, if both $X$ and $Y$ are $B$-submodule coalgebras of $C$, then so is $X \land Y$. \\

For any  collection $\RR$ of $B$-module coalgebras, we simply call a member of $\RR$ an
$\RR$-{\em coalgebra}. For $C \in \Coa{B}$, a
$B$-submodule coalgebra of $C$, which is also a member of $\RR$, is called an {\em
$\RR$-subcoalgebra} of $C$. Now, we can introduce our
definition of a radical class.

\begin{defn}
Let $B$ be a bialgebra. A non-empty collection $\RR$ of objects in $\Coa{B}$
is said to be a {\em radical class} if it satisfies the following
conditions:
\begin{enumerate}
\item[(R1)] $\RR$ is closed under subobjects, i.e. every
$B$-submodule coalgebra of an $\RR$-coalgebra is an $\RR$-coalgebra.
\item[(R2)] For any $C \in \Coa{B}$ and $\{C_i\}_{i \in I}$ a family
of $\RR$-subcoalgebras of $C$,  $\sum_{i \in I} C_i \in \RR$.
\end{enumerate}
A radical class $\RR$ of $\Coa{B}$ is called {\em idempotent} if it
satisfies
\begin{enumerate}
\item[(R3)]  For any $C \in \Coa{B}$,  if $C_1$, $C_2$ are $\RR$-subcoalgebras of $C$,
then $C_1 \land C_2 \in\RR$.
\end{enumerate}
\end{defn}
\begin{defn}
  Let $\RR$ be a radical class of $\Coa{B}$. For $C \in \Coa{B}$, the
  $\RR$-\emph{radical} of $C$ is defined as the sum of all the $\RR$-subcoalgebras of $C$, and
  is denoted
  by $\RR(C)$.
\end{defn}

It follows immediately from the definition that $\RR(C)$ is the
unique maximal $\RR$-subcoalgebra of the $B$-module coalgebra $C$.
Indeed, the notion of $\RR$-radical of module coalgebras is a
generalization of coradical of coalgebras.

\begin{example}
Let $B$ be the trivial bialgebra $k$. Then $\Coa{B}$ is simply the
category of all coalgebras over $k$. Consider the collection $\RR$ of all
semisimple coalgebras over $k$. By \cite[Theorem 3]{Kapl75}, every
subcoalgebra of a semisimple coalgebra is semisimple, and  the sum
of a collection of semisimple subcoalgebras of a coalgebra $C$ is also semisimple.
Therefore, $\RR$ is a radical class of $\Coa{B}$, and $\RR(C)$ is
the coradical $C_0$ of $C$. Obviously, $\RR$ is not idempotent.
\end{example}
Similar to the properties of coradical of coalgebras, we have the following result:
\begin{prop}
 Let $\RR$ be a radical class of $\Coa{B}$, and $C \in \Coa{B}$.
\begin{enumerate}
\item[\rm (i)] If $D$ is a $B$-submodule coalgebra of $C$, then
 $\RR(D) = \RR(C) \cap D$.
\item[\rm (ii)] If $C$ is a direct sum of a family $\{C_i\}$ of $B$-submodule coalgebras,
then $\RR(C) = \bigoplus_i\RR(C_i)$.
\end{enumerate}
\end{prop}
\begin{proof}
(i) Clearly, $\RR(D) \subseteq \RR(C)$ and hence $\RR(D) \subseteq
\RR(C) \cap D$. Conversely, since $\RR(C) \cap D$ is an $\RR$-subcoalgebra of $D$, $\RR(C) \cap D \subseteq \RR(D)$.

(ii)  If $C=\bigoplus_i C_i$ as $B$-module
coalgebras, it follows from \cite[Theorem 3]{Kapl75} that $\RR(C)
= \bigoplus_i (\RR(C) \cap C_i)$. Hence, by (i), $\RR(C) =
\bigoplus_i \RR(C_i)$.
\end{proof}

If $\RR$ is an idempotent radical class of $\Coa{B}$, it is immediate to see that
$$
\RR(C) \wedge \RR(C) = \RR(C)\,.
$$
Moreover, the membership of a $B$-module coalgebra in $\RR$ can be
determined by its simple subcoalgebras.

\begin{prop}\label{c-radical2}
Let $\RR$ be an idempotent radical class of $\Coa{B}$.
 The following statements concerning a $B$-module coalgebra $C$ are equivalent :
\begin{enumerate}
\item[\rm (i)] $C$ is an $\RR$-coalgebra.
\item[\rm (ii)] There exists an $\RR$-subcoalgebra $C'$ of $C$ which contains
the coradical $C_0$ of $C$.
\item[\rm (iii)]  For every simple
subcoalgebra $D$ of $C$, $DB$ is an  $\RR$-subcoalgebra  of $C$.
\end{enumerate}
\end{prop}
\begin{proof}
 ((i) $\Rightarrow$ (ii)) and   ((ii) $\Rightarrow$ (iii)) follow directly from (R1). \\
(iii) $\Rightarrow$ (i): Since $DB$ is an $\RR$-subcoalgebra  of $C$
for every simple subcoalgebra $D$ of $C$, $C_0B$ is an
$\RR$-subcoalgebra by property (R2) where $C_0$ is the coradical of
$C$. It follows from (R3) that $\bigwedge^n C_0B$ is an
$\RR$-subcoalgebra for any $n \ge 1$. Hence, $\sum_{n\ge 1}
\bigwedge^n C_0B$ is an $\RR$-subcoalgebra. The result follows from (R2) and
the fact that
$$
C=\sum_{n\ge 1} \sideset{}{^n}\bigwedge C_0B. \qedhere
$$
\end{proof}

In closing this section, we give several equivalent formulations of
idempotent radical class for the sake of convenience of further discussion.

\begin{prop}\label{p:e3}
Let $B$ be a bialgebra, and $\RR$ a non-empty collection of objects in
$\Coa{B}$ which satisfies the conditions
{\rm (R1)} and {\rm (R3)}. Then $\RR$ is an idempotent radical class
of $\Coa{B}$ if, and only if, one the following conditions holds:
\begin{enumerate}
  \item[(R2)] For any $C \in \Coa{B}$,  if $\{C_i\}_{i \in I}$ is a family of $\RR$-subcoalgebras of $C$, then
  $\sum_i C_i$ is an $\RR$-subcoalgebra of $C$.
  \item[(R2a)] For any $C \in \Coa{B}$,  if $\{C_i\}_{i \in I}$ is a chain of
   $\RR$-subcoalgebras of $C$, then $\bigcup_i C_i$ is also an $\RR$-subcoalgebra of $C$.
  \item[(R2b)] For any $C \in \Coa{B}$, there exists a maximal $\RR$-subcoalgebra of $C$.
\end{enumerate}
\end{prop}
\begin{proof} It suffices to show that (R2), (R2a) and (R2b) are equivalent conditions under the hypotheses (R1)
and (R3). Obviously, (R2a) is an immediate consequence of (R2). By Zorn's Lemma, (R2b) follows from (R2a).
  Suppose that the condition (R2b) holds for $\RR$. Let $C \in \Coa{B}$ and $D$ a maximal $\RR$-subcoalgebra of $C$.
  We claim that $D$ contains all the $\RR$-subcoalgebras of $C$. If not, there exists an $\RR$-subcoalgebra
  $D'$ of $C$ such that $D' \nsubseteq D$. Then
  $$
  D \subsetneq D' + D \subseteq D'\land D.
  $$
  By (R3),
  $D' \land D$ is also an $\RR$-subcoalgebra of $C$ which properly contains $D$.
  This contradicts the maximality of $D$.

  Let $\{C_i\}$ be a family of $\RR$-subcoalgebras of $C$. Then $C_i \subseteq D$ for $i \in I$. Therefore,
  $\sum_{i \in I}C_i \subseteq D$. It follows from (R1) that $\sum_{i \in I}C_i$ is also an $\RR$-subcoalgebra of
  $C$. Hence,  (R2) holds for $\RR$.
\end{proof}

\section{Projective Radical}
In this section, we will show that if $H$ is a Hopf algebra, the collection $\PP$
of ($H$-)projective $H$-module coalgebras
is an idempotent radical class of $\Coa{H}$. We begin with a short review on Hopf modules. \\

Let $B$  be a bialgebra and $C \in \Coa{B}$. Recall that a (right)
$(C, B)$-Hopf module $M$ is a right $B$-module and a right
$C$-comodule such that the  $C$-comodule structure map $\rho: M
\to M \otimes C$ is a $B$-module map, where $M \otimes C$ uses the diagonal $B$-action.
The category of  $(C,B)$-Hopf
modules is denoted by $\C{B}^C$. The morphisms between two
$(C,B)$-Hopf modules are those $B$-module maps which are also
$C$-comodule maps at the same time. The category ${^C\!\C{B}}$ can also be defined
similarly, and both  ${^C\!\C{B}}$ and $\C{B}^C$  are
abelian categories. \\

Recall the ``wedge'' product for comodules from \cite{Sw69}. Suppose
$C \in \Coa{B}$ and $M \in \M_B^C$. For $k$-subspaces $N \subseteq M$ and  $X\subseteq C$,
the wedge product $N\land X$ is defined to be the
kernel of the map
$$
M \stackrel{\rho_M}{\to} M \otimes C  \to  M/N \otimes C/X\,.
$$
If $N$ is a
$B$-submodule of $M$ and $X$ is a
$(C,B)$-Hopf submodule of $C$, then
 $N \land X$ is also  a $(C,B)$-Hopf submodule of $M$. If $M=C$ and both
$N$ and $X$ are $B$-submodule coalgebras of $C$, then $N\land
X$ is a $B$-submodule coalgebra of $C$.
\begin{lem}\label{r:l1}
Let $B$ be a bialgebra and $C \in \Coa{B}$. Suppose that $C_1$, $C_2$ are $B$-submodule coalgebras of $C$.
Then for any $M \in \M_B^{C_1 \land C_2}$,
the sequence
$$
0 \to \{0\}\land C_1 \to M \to M/(\{0\}\land C_1) \to 0
$$
is exact in $\M^{C_1 \land C_2}_B$. Moreover,
 $\{0\} \land C_1 \in \M^{C_1}_B$ and $M/(\{0\}\land C_1) \in
\M^{C_2}_B$.
\end{lem}
\begin{proof}
Notice that $C_1$, $C_2$ are $B$-submodule coalgebras of $C_1
\land C_2$.  Hence, $\{0\}\land C_1$ is a $(C_1\land C_2, B)$-Hopf submodule of $M$.
Since $\M^{C_1 \land C_2}_B$ is an abelian category, the quotient
$M/(\{0\}\land C_1)$ is also a $(C_1\land C_2, B)$-Hopf module and the exact sequence
follows easily. It follows immediately from the definition of wedge product that
$$
\rho_M(\{0\}\land C_1) \subseteq (\{0\}\land C_1)  \otimes C_1 \,,
$$
where $\rho_M: M \to M \otimes C$ is the $(C_1 \land C_2)$-comodule structure map of $M$.
Therefore, $\{0\}\land C_1 \in \M^{C_1}_B$. Since
$$
M =  \{0\} \land (C_1 \land C_2)\,,
$$
 by the associativity of wedge product (cf. \cite[Proposition 9.0.0]{Sw69}),
$$
M = (\{0\} \land C_1) \land C_2 \,.
$$
Let $\eta_1 : M \to M/(\{0\} \land C_1)$ and $\eta_2 : C_1 \land C_2 \to (C_1 \land C_2)/C_2$
be the natural surjections. Then
$$
(\eta_1 \otimes \eta_2)\circ\rho_M = 0\,.
$$
Hence,
$$
(\eta_1 \otimes 1)\circ\rho_M(M) \subseteq M/(\{0\}\land C_1) \otimes C_2 \,.
$$
Therefore, $M/(\{0\}\land C_1) \in \M^{C_2}_B$.
\end{proof}

\begin{thm}\label{t:pr}
Let $H$ be a Hopf algebra, and $\PP$ the collection of all right $H$-module coalgebras which are
  projective  $H$-modules. Then $\PP$ is an idempotent radical class of $\Coa{H}$.
\end{thm}
\begin{proof}
Let $C \in \Coa{H}$.
          Recall from \cite{Doi83} that $C$ is a projective $H$-module if, and only if, there
          exists a right $H$-module map $\psi: C \to H$ such that
          \begin{equation}\label{*}
            \e_H \circ \psi  = \e_C.
          \end{equation}
          In this case
          all the  right $(C, H)$-Hopf modules are projective $H$-modules. By Proposition \ref{p:e3}, it suffices
          to show that $\PP$ satisfies (R1), (R3) and (R2b).

          (R1) Let $C$ be a $\PP$-coalgebra. For any $H$-submodule coalgebra $D$ of $C$, $D$ is a right
          $(C,H)$-Hopf module and hence a projective $H$-module.

          (R3) Let $C \in \Coa{H}$, and let $C_1$, $C_2$ be $\PP$-subcoalgebras of
          $C$. Then, by the preceding remark, $C_1 \land C_2$ is also a right $H$-module coalgebra.
          Now, by Lemma \ref{r:l1}, the
            sequence
            $$
                0 \to C_1 \to C_1 \land C_2 \to (C_1 \land C_2)/C_1 \to 0
            $$
            is exact in $\C{H}^{C_1 \land C_2}$. Moreover, $(C_1 \land C_2)/C_1 \in \C{H}^{C_2}$.
            Therefore, $(C_1 \land C_2)/C_1$ is a projective $H$-module, and so
            the sequence is split exact in $\CMod{H}$. Thus $C_1 \land C_2 \cong  C_1 \oplus (C_1 \land C_2)/C_1$
            as $H$-modules. Hence $C_1\land C_2$ is a projective $H$-module.

           (R2b) Let $C \in \Coa{H}$ and $S$ the set of all pairs $(D, \psi)$ in which $D$ is
           a $\PP$-subcoalgebra of $C$ and $\psi: D \to H$ is a right $H$-module map satisfying \eqref{*}.
           We define the partial ordering $\le$ on the non-empty set $S$ as follows:
           $$
           (D, \psi) \le (D', \psi')\quad \text{if}\quad D \subseteq D'\quad \text{and} \quad \psi'|_D = \psi\,.
           $$
           Suppose that $\{(D_i, \psi_i)\}$ is a chain in $S$. Let $\ol D=\bigcup_i D_i$ and define $\ol \psi : \ol D \to H$ as
           $$
           \ol \psi(x) = \psi_i(x) \quad \text{if}\quad x \in D_i\,.
           $$
           Obviously, $\ol D$ is a right $H$-submodule coalgebra of $C$.
           Since $\{(D_i, \psi_i)\}$ is a chain, the function
           $\ol\psi: \ol D \to H$ is a well-defined right $H$-module map which satisfies \eqref{*}, and so
           $(\ol D, \ol \psi) \in S$.
           By
           Zorn's Lemma, there is a maximal element $(D, \psi) \in S$.
           If follows from \cite{Doi83} that $D$ is a $\PP$-subcoalgebra of $C$.
           We now claim that $D$ is a maximal
           $\PP$-subcoalgebra of $C$. If the claim is false, then there exists a $\PP$-subcoalgebra $D'$ of $C$ such that
           $D \subsetneq D'$. Since $D'/D \in \C{H}^{D'}$, $D'/D$ is a right projective $H$-module and so the
           sequence
           $$
           0 \to D \to D' \to D'/D \to 0
           $$
           is split exact in $\CMod{H}$. Therefore, $D'=D \oplus M$ for some $H$-submodule $M$ of $D'$. Let
           $\psi': D' \to H$ be a right $H$-module map satisfying \eqref{*}. Consider the map
           $\hat{\psi}=\psi \oplus \psi'|_M$. Then $\hat{\psi}: D' \to H$ is also a right $H$-module map
           satisfying \eqref{*} and $\hat{\psi}|_D=\psi$. This leads to $(D, \psi) \lvertneqq (D', \hat{\psi})$,
           a contradiction!
\end{proof}

The following corollary generalizes \cite[Corollary 10]{Ng98} to arbitrary Hopf algebras.

\begin{cor}\label{c:p}
Let $H$ be a Hopf algebra. For every   $H$-module coalgebra $C$,
there exists a unique maximal projective $H$-submodule
coalgebra $\PP(C)$ of $C$.
Moreover, the following statements about a right $H$-module coalgebra $C$ are equivalent:
\begin{enumerate}
  \item[\rm (i)] $C$ is $H$-projective;
  \item[\rm (ii)] there exists an $H$-projective submodule coalgebra of $C$ which contains the
  coradical $C_0$;
  \item[\rm (iii)] $DH$ is
$H$-projective for every simple subcoalgebra $D$ of $C$.
\end{enumerate}
\end{cor}
\begin{proof}
  The corollary follows directly from Theorem \ref{t:pr} and Proposition \ref{c-radical2}.
\end{proof}

As an application of projective radical of $H$-module coalgebra, we give another proof for a
projectivity result obtained by Schneider in
\cite{Schn90}.

\begin{cor}[Schneider]
Let $H$ be a Hopf algebra and $C$ a right $H$-module coalgebra over
the field $k$. Let $G(C)$ be the set of group-like elements of $C$. Assume that
\begin{enumerate}
\item[\rm (i)] the coradical of $C \otimes k'$ is contained in  $G(C \otimes k')(H \otimes k')$ for an extension $k'
\supseteq k$, and
\item[\rm (ii)] the canonical map $can : C \otimes H \to C \otimes
C$ defined by
$$
 can : c \otimes b \mapsto \sum c_1 \otimes c_2 b
$$
is injective.
\end{enumerate}
Then $C$ is a projective $H$-module.
\end{cor}
\begin{proof}
Since $C$ is $H$-projective if, and only if, $C \otimes k'$ is $H \otimes
k'$-projective for some extension $k' \supseteq k$, we may simply
assume $k=k'$.

For $g \in G(C)$, the surjective
map $\phi : H \to gH$ defined by $\phi(h) = gh$ is an $H$-module
coalgebra map.  The injectivity of $can$ implies that $\phi$ is also
an isomorphism. In particular, $gH$ is an $H$-free submodule coalgebra of $C$.
By Theorem \ref{t:pr} and (R2), $G(C)H=\sum_{g\in G(C)}gH$ is an $H$-projective submodule coalgebra of $C$,
and it contains the coradical of $C$.  It follows from Corollary
\ref{c:p} that $C$ is a projective $H$-module.
\end{proof}
\section{$\EE$ is an idempotent Radical Class}
Let $C$ be a $B$-module coalgebra. Then
$\overline{C}= C/CB^+$ admits a natural coalgebra structure, where $B^+= \ker \e_B$, and the
natural surjection $\eta_C: C \to \overline{C}$ is a coalgebra map.
For any right $B$-module $M$, we let
$$
\overline{M} =  M/(MB^+),
$$
and $\eta_M :M \to \overline{M}$, $m \mapsto
\overline{m}$, the quotient map. If $M \in \M_B^C$, then $\overline{M}$ admits a
natural  right $\Cbar$-comodule structure, and $\ol ?: \C{B}^C \to \M^{\ol C}$, $M \mapsto \ol M$,
defines a $k$-linear functor.\\

 For any $N \in \M^\Cbar$,
the cotensor product $N \square_\Cbar C$ is defined as the kernel of the
map
$$
\rho_N \otimes \id_C - \id_N \otimes (\eta_C \otimes \id)\circ\Delta_C : N \otimes C \to N \otimes {\ol C} \otimes C.
$$
The cotensor product $N \square_{\ol C} C$ has a right $(C,B)$-Hopf module
structure inherited from $C$, and
$?\,\square_\Cbar C : \M^{\ol C} \to \C{B}^C$ is a $k$-linear functor. Moreover,
the functor $\ol{?}$ is left
adjoint to $?\,\square_\Cbar C$ with the unit $\Xi$ and counit
$\Theta$ of the adjunction  given by
\begin{eqnarray*}
\Xi_M : M \to \overline{M}\square_\Cbar C &,& m \mapsto \sum
\overline{m_0} \otimes m_1 \\
\Theta_N : \overline{N \square_\Cbar C} \to N &,& \overline{\sum n_i
\otimes c_i} \mapsto \sum n_i \e(c_i)
\end{eqnarray*}
where $M \in \M_B^C$ and $N \in \M^\Cbar$. In particular,
$?\,\square_{\ol C} C$ is left exact and $\ol ?$ is right exact (cf. \cite[Theorem 2.6.1]{Weibhobk}). \\

In this section, we prove in Corollary \ref{c:e3} that $\ol ?: \C{B}^C \to \M^{\ol C}$ is an equivalence
if, and only if, for all $M \in \C{B}^C$, (i) $\Tor(M,k)=0$,
and (ii) $M=0$ whenever $\ol M =0$.
Using this characterization for the equivalence of $\ol ?$, we show in
Theorem \ref{t:eclass} that the collection
$$
\EE=\{ C \in \Coa{B} \mid \text{ the functor } \ol ?: \C{B}^C \to \M^{\ol C} \text{ is an equivalence}\}
$$
is an idempotent radical class of $\Coa{B}$. The result will be used to show that the $B$-cocleft
module coalgebras also form an idempotent radical class. We begin with a characterization for the
equivalence of the functor $\ol ?$.

\begin{prop}\label{p:e1}
Let $B$ be a bialgebra, and $C$ a right $B$-module coalgebra. Then the
following statements are equivalent:
\begin{enumerate}
  \item[(i)] the functor $\ol ? : \C{B}^C \to \M^{\ol C}$ is an
  equivalence;
  \item[(ii)] the functor $\ol ?$ is exact, and for any $M \in \C{B}^C$, $\ol M =0$ implies $M=0$;
  \item[(iii)] the unit $\Xi$ and counit $\Theta$ of the adjunction  are isomorphisms.
\end{enumerate}
\end{prop}
\begin{proof}
The implications ((i) $\Rightarrow$ (ii)) and ((iii) $\Rightarrow$
(i))
are straightforward. It remains to show that (ii) implies (iii).

For any $k$-space $V$, $V \otimes C$ admits a right $(C,B)$-Hopf module structure induced by the $B$-action and
comultiplication of $C$, and $\FF (V):=V \otimes
\ol C$ is a right $\ol C$-comodule with its structure inherited from the
comultiplication of $\ol C$. The map
$$
\FF(V) \square_{\ol C} C \xrightarrow{\id_V \otimes \e \otimes \id_C} V \otimes C
$$
is an isomorphism of $(C,B)$-Hopf modules. Note that
$\tilde\eta: \ol{V \otimes C} \to \FF(V)$, defined by $\tilde\eta: \ol{v \otimes c}\mapsto v \otimes \ol c$,
is an isomorphism of right $\ol C$-comodules, and the diagram
$$
\xymatrix{
\ol{\FF(V) \square_{\ol C} C} \ar[rd]^-{\Theta_{\FF(V)}}\ar[rr]^-{\ol{\id_V \otimes \e \otimes \id_C}}
&& \ol{V \otimes C} \ar[ld]_-{\tilde\eta} \\
&\FF(V)&
}
$$
commutes. Therefore, $\Theta_{\FF (V)}$ is an isomorphism.

For $N \in
\M^{\ol C}$, the  $\ol C$-comodule structure map $\rho_N : N
\to \FF(N)$ is an injective $\ol C$-comodule map. Let
$$
N' := \coker \rho_N = \FF(N)/\rho_N(N),
$$
and let $g$ be the composition
$$
\FF(N) \xrightarrow{\pi} N' \xrightarrow{\rho_{N'}} \FF(N'),
$$
where $\pi$ is the natural surjection map. Then we have the exact
sequence
$$
0 \to N \xrightarrow{\rho_N} \FF(N) \xrightarrow{g}  \FF(N')
$$
in $\M^{\ol C}$.

Assume (ii) holds. Since $?\,\square_{\ol C} C$ is left exact and
$\ol ?$ is exact, we have the following commutative diagram:
$$
\xymatrix{ 0\ar[r] & \overline{N \square_\Cbar C} \ar[r]
\ar[d]^{\Theta_N} & \overline{\FF(N)\square_\Cbar C} \ar[r]
\ar[d]^-{\Theta_{\FF(N)}} & \overline{\FF(N')\square_\Cbar C}\ar[d]^-{\Theta_{\FF(N')}} \\
0\ar[r] & N  \ar[r]  & \FF(N) \ar[r] & \FF(N')
}
$$
where the rows are exact, and both  $\Theta_{\FF(N)}$ and
$\Theta_{\FF(N')}$ are isomorphisms. It follows from diagram tracing that $\Theta_N$
is also an isomorphism.

For any $M \in \M^C_B$, we have the exact sequence
$$
\xymatrix{ 0 \ar[r] &\ker \Xi_M \ar[r]& M \ar[r]^-{\Xi_M} &
\overline{M} \square_{\Cbar} C \ar[r] & \coker \Xi_M \ar[r] & 0}
$$
in $\M_B^C$. By the exactness of $\ol{?}$, we also have the exact
sequence
$$
\xymatrix{ 0 \ar[r] & \overline{\ker \Xi_M}\ar[r] & \overline{M}
\ar[r]^-{\ol {\Xi}_M} & \overline{\overline{M}\square_{\Cbar} C}
\ar[r] & \overline{\coker \Xi_M}\ar[r] & 0 }
$$
in $\M^{\Cbar}$. Since the equation
$$
\Theta_{\ol M} \circ  \ol {\Xi}_M = \id_{\ol M}
$$
holds for every $M \in \C{B}^C$ and $\Theta_{\ol M}$ is an
isomorphism, so is $\ol \Xi_M$. Therefore,
$$
\overline{\ker \Xi_M} = 0 = \overline{\coker \Xi_M}\,.
$$
It follows from (ii) that
$$
\ker \Xi_M = 0 = \coker \Xi_M.
$$
Hence, $\Xi_M$ is an isomorphism.
\end{proof}
\begin{prop}\label{p:e2}
  Let $B$ be a bialgebra and $C$ a right $B$-module coalgebra. Then
  the functor $\ol ? : \C{B}^C \to \M^{\ol C}$ is exact if, and only
  if $\Tor(M, k)=0$ for all $M \in \C{B}^C$, where $k$ is considered
  as a trivial left $B$-module.
\end{prop}
\begin{proof}
Let
$$
0 \to M_1 \xrightarrow{a} M_2 \xrightarrow{b} M_3 \to 0
$$
be an exact sequence in $\C{B}^C$. By the associated long exact sequence, the sequence
\begin{equation}\label{eq:tor-exact}
\Tor(M_3, k) \xrightarrow{} M_1 \otimes_B k \xrightarrow{a \otimes_B k}
M_2\otimes_B k \xrightarrow{b \otimes_B k}  M_3 \otimes_B k \xrightarrow{} 0
\end{equation}
is exact. Note that $\phi_M: \ol M \to M \otimes_B k$, defined by $\phi: \ol m \mapsto
m \otimes 1$ for $m \in M$, is a natural isomorphism of $k$-linear spaces.
If $\Tor(M, k)=0$ for all $M \in \C{B}^C$, then it follows from \eqref{eq:tor-exact} that
the sequence
$$
0 \xrightarrow{} M_1 \otimes_B k \xrightarrow{a \otimes_B k}
M_2\otimes_B k \xrightarrow{b \otimes_B k}  M_3 \otimes_B k \xrightarrow{} 0
$$
is exact. By the naturality of $\phi$ that the sequence
$$
0\to \ol{M_1} \xrightarrow{\ol a}  \ol{M_2}\xrightarrow{\ol b}\ol{M_3} \to 0
$$
is exact in $\M^{\ol C}$. Therefore, $\ol ?$ is an exact functor.

Conversely, assume $\ol ?$ is exact. For $M \in \C{B}^C$ with $C$-comodule structure $\rho_M: M \to M\otimes C$,
let $\GG(M):=M \otimes B$ be the right $(C,B)$-Hopf module with
the $B$-action $\cdot$, and the $C$-coaction $\rho_{\GG(M)} : \GG(M) \to \GG(M) \otimes C$ given by
$$
\begin{aligned}
  (m \otimes b)\cdot h & = m \otimes bh\quad \text{and} \\
  \rho_{\GG(M)}(m \otimes b) &= \sum m_0 \otimes b_1 \otimes m_1b_2
\end{aligned}
$$
for all $m \in M$ and $b, h \in B$, where $\rho_M(m)=\sum m_0 \otimes m_1$ and $\Delta(b) =\sum b_1 \otimes b_2$.
Notice that the $B$-module structure map  $\mu: \GG(M)  \to M$ of $M$ is a $(C,B)$-Hopf module map, and so
we have
the exact sequence
$$
 \xymatrix{ 0 \ar[r] & \ker \mu \ar[r]^-{i} & \GG(M) \ar[r]^-{\mu} & M \ar[r] & 0}
$$
in $\C{B}^C$, where $i$ is the inclusion map.
By the naturality of $\phi$, we have the commutative diagram
\begin{equation}\label{eq:cd1}
 \xymatrix{ & 0 \ar[r] &
\ol{\ker \mu} \ar[r]^-{\ol i} \ar[d]^-{\phi} & \ol{\GG(M)} \ar[r]^-{\ol \mu}
\ar[d]^-{\phi}
& \ol{M}  \ar[d]^-{\phi}  \\
\Tor(\GG(M), k) \ar[r]&\Tor(M, k) \ar[r] & \ker \mu \otimes_B k \ar[r]^-{i \otimes_B k} &
\GG(M)\otimes_B k \ar[r]^-{\mu \otimes_B k} & M \otimes_B k \,,}
\end{equation}
where top and bottom rows are exact by the exactness of $\ol ?$ and
the associated long exact sequence respectively. In particular, $\ker (i \otimes_B k) =0$ and so the map
$$
\Tor(\GG(M), k) \to \Tor(M, k)
$$
in the diagram is surjective. Since $\GG(M)$ is a free right $B$-module,
$\Tor(\GG(M),k)=0$.     Therefore, $\Tor(M,k)=0$.
\end{proof}
\begin{cor}\label{c:e3}
  Let $B$ be a bialgebra, and $C$ a  $B$-module coalgebra. Then the functor $\ol ?: \C{B}^C \to \M^{\ol C}$ is an equivalence
  if, and only if, the following conditions hold for all $M \in \C{B}^C$:
  \begin{enumerate}
    \item[(i)] $\Tor(M,k)=0$, and
    \item[(ii)] $\ol M=0$ implies $M=0$.
  \end{enumerate}
\end{cor}
\begin{proof}
  The statement follows immediately from Propositions \ref{p:e1} and \ref{p:e2}.
\end{proof}
Now we turn to the main result of this section.
\begin{thm}\label{t:eclass}
  Let $B$ be a bialgebra and $\EE$ the collection of all the  $B$-module coalgebras $C$
  such that $\ol ? : \C{B}^C \to \M^{\ol C}$ is an equivalence. Then $\EE$ is an idempotent radical class of
  $\Coa{B}$.
\end{thm}
\begin{proof} By Proposition \ref{p:e3}, it suffices to show that $\EE$  satisfies (R1), (R2a) and (R3).

  (R1) Let $C$ be an $\EE$-coalgebra, and $D$ a $B$-submodule coalgebra of $C$. For $M \in \C{B}^D$,
  $M \in \C{B}^C$. By Corollary \ref{c:e3}, $\ol? : \C{B}^D \to \M^{\ol D}$ is an equivalence.

  (R2a) Let $C$ be a right $B$-module coalgebra and $\{C_i\}$ a chain of $\EE$-subcoalgebras of $C$.
     Then $D:=\bigcup_i C_i$ is a $B$-submodule coalgebra of $C$. For $M \in \C{B}^D$, let
     $M_i=\{0\} \land C_i$. Then $M_i \in \C{B}^{C_i}$ and
    \begin{equation}\label{eq:union}
      M =\bigcup_i M_i\,.
    \end{equation}
    It follows from Corollary \ref{c:e3} that
    $$
    \Tor(M_i, k)=0
    $$
    for all $i$. Since $\Tor(-, k)$ commutes with direct limit (cf. \cite[Corollary 2.6.17]{Weibhobk}), we have
    $$
    \Tor(M,k) = \lim_{\rightarrow}\, \Tor(M_i, k) = 0\,.
    $$
    Therefore, the functor $\ol ?: \C{B}^D \to \M^{\ol D}$ is exact by Proposition \ref{p:e2}.
    Suppose $\ol M = 0$. By the exactness of $\ol ?$, the sequence
    $$
    \xymatrix{
    0 \ar[r] & \ol M_i \ar[r]^-{\ol j} & \ol M
    }
    $$
    is exact for all $i$ , where $j$ is the inclusion map. Hence, $\ol M_i =0$ for all $i$. By
    Corollary \ref{c:e3}, $M_i = 0$ for all $i$. It follows from \eqref{eq:union} that $M=0$.
    Therefore, by Corollary \ref{c:e3}, $D$ is an $\EE$-subcoalgebra of $C$.

    (R3) Let $C_1$, $C_2$ be $\EE$-subcoalgebra of a $B$-module coalgebra $C$, and
    $D'=C_1 \land C_2$. For $M \in \C{B}^{D'}$, by Lemma \ref{r:l1}, we have
    the exact sequence
    \begin{equation}\label{eq:filteration}
     0 \to \{0\}\land C_1 \to M \to M/(\{0\}\land C_1) \to 0
    \end{equation}
    in $\M^{D'}_B$ with $\{0\} \land C_1 \in \M^{C_1}_B$ and $M/(\{0\}\land C_1) \in
     \M^{C_2}_B$. In particular,
     $$
     \Tor(M/(\{0\}\land C_1), k) = 0 =\Tor(\{0\}\land C_1, k)\,.
     $$
     By the associated long exact sequence of \eqref{eq:filteration}, we have the exact sequence
     $$
     \Tor(\{0\}\land C_1, k) \to \Tor(M, k) \to \Tor(M/(\{0\}\land C_1), k)\,.
     $$
     Therefore, $\Tor(M,k)=0$. Hence by Proposition \ref{p:e2},
     the functor $\ol ?:\C{B}^{D'} \to \M^{\ol {D'}}$ is exact. Suppose $\ol M=0$. Then, by \eqref{eq:filteration}
     and the exactness of $\ol ?$, we have
     $$
     \ol{\{0\} \land C_1} = 0 =\ol{M/(\{0\}\land C_1)}.
     $$
     Since $C_1$, $C_2$ are $\EE$-coalgebras, it follows from Corollary \ref{c:e3} that
     $$\{0\} \land C_1 =0= M/(\{0\}\land C_1)$$
     and hence $M=0$. By  Corollary \ref{c:e3} again, $D'$ is an $\EE$-coalgebra.
\end{proof}
\section{Cocleft Radical}
In this section, we prove that the collection $\CC$ of all cocleft $B$-module coalgebras is an idempotent radical class.
As an application, we use this result to give another proof for a normal basis theorem of Schneider
\cite[Theorem III]{Schn90} without assuming bijective antipode of the underlying Hopf algebra. We begin with
some definitions and properties of cocleft $B$-module coalgebras.
\begin{defn}
Let $B$ be a bialgebra and $C$ a right $B$-module coalgebra.
\begin{enumerate}
\item[(i)] $C$ is said to have a {\em normal basis} if
$$
C \cong \overline{C} \otimes B \quad \text{as left $\ol C$-comodule
and right $B$-module.}
$$
\item[(ii)] $C$ is said to be {\em $B$-cogalois} if the  map
$$
can : C \otimes B \to C {\square_\Cbar} C \,,\quad  c \otimes b
\mapsto \sum c_1 \otimes c_2 b
$$
is a $k$-linear isomorphism.
\item[(iii)] A $B$-module map $\gamma :C\to B$ is called a {\em
cointegral} if $\gamma$ is convolution invertible.
\item[(iv)] $C$ is said to be {\em $B$-cocleft} it admits a
cointegral.
\item[(v)] A $B$-submodule coalgebra $D$ of $C$ is called
  a $B$-cocleft subcoalgebra if $D$ is $B$-cocleft.
\end{enumerate}
\end{defn}

It was proved in \cite{MaDo92} that $C$ is $B$-cocleft if, and only
if, $C$ is $B$-cogalois and has a normal basis. In this case, $\ol ?: \M_B^C
\to \M^\Cbar$ is an equivalence and all the $(C, B)$-Hopf modules are free
$B$-modules. In particular, a $B$-cocleft coalgebra is an $\EE$-coalgebra or $\CC \subset \EE$.

\begin{lem}\label{l:4.2}
Let $B$ be a bialgebra, and $C$ an $\EE$-coalgebra.
Suppose that for every simple subcoalgebra $E$ of
$\ol{C}$ there exists a $B$-cocleft subcoalgebra $D$ of $C$
such that $E \subseteq \eta_C(D)$, where $\eta_C: C \to \ol C$ is the natural surjection.
Then $C$ is $B$-cocleft.
\end{lem}
\begin{proof}
By Proposition \ref{p:e1}, both of $\Theta$ and $\Xi$ are isomorphisms and so $?\,\square_{\ol{C}} C
:\M^{\ol{C}} \to \M_B^C$ is  a $k$-linear  equivalence. In
particular, $C$ is an injective left $\ol{C}$-comodule. In order to
show that $C$ is $B$-cocleft, by \cite[Theorem 2.3]{MaDo92}, it
suffices to show that $C$ admits a normal basis. By \cite[Collorary 2.2]{Schn90}, it is enough to
show that $E \square_{\ol C} C \cong E \otimes B$ as left $E$-comodule and right $B$-module for all
simple subcoalgebra $E$ of $\ol C$.

Let $E$ be a simple
subcoalgebra of $\ol{C}$, and $D$ a $B$-cocleft subcoalgebra of $C$
such that $E \subseteq \eta_C(D)$. Let $i: D \to C$ be the inclusion map. Then $i$ is a monomorphism in $\C{B}^C$ and
${\ol i} ({\ol D})= \eta_C(D)$.
Since $\ol ?: \C{B}^C \to \M^{\ol C}$ is exact, $\ol i: \ol D \to \eta_C(D)$ is a coalgebra
isomorphism. Thus $E$ can be viewed as a right $\ol D$-comodule
via $\ol i$.
 Notice that $E \square_{\ol{D}} D$ is a $(C, B)$-Hopf submodule of
 $E \square_{\ol{C}} C$. Let $j: E \square_{\ol{D}} D \to E \square_{\ol{C}}
 C$
 be the inclusion map. Then we have the following commutative diagram
 $$
\xymatrix{ \ol{E \square_{\ol{D}} D} \ar[r]^-{\ol{j}}
\ar[d]_-{\Theta_E} & \ol{E \square_{\ol{C}} C}
\ar[d]^-{\Theta_E}\\
 E \ar@{=}[r]& E
 }\,
$$
 of $k$-linear maps. Since both $D$ and $C$ are $\EE$-coalgebra, the vertical maps $\Theta_E$ in the
 diagram are isomorphisms by Proposition \ref{p:e1}.
  Hence, $\ol j$ is an isomorphism. Since $\ol ?$ is an equivalence, $j$ is
 also  an isomorphism and so
$$
E \square_{\overline{C}} C  =  E \square_{\overline{D}} D \,.
$$
Since $D$ has a normal basis, it follows from  \cite[Corollary
2.2]{Schn90} that
$$
E \square_{\overline{C}} C = E \square_{\overline{D}} D \cong E \otimes B
$$
as left $E$-comodule and right $B$-module.
\end{proof}
Now we can prove our main theorem of this section.
\begin{thm}\label{t:cocleft}
  Let $B$ be a bialgebra. Then the collection $\CC$ of all
  $B$-cocleft coalgebras is an idempotent radical
  class of $\Coa{B}$.
\end{thm}
\begin{proof} By Proposition \ref{p:e3}, we only need to show that $\CC$ satisfies (R1), (R2a) and (R3).

  (R1) Let $C$ be a $B$-cocleft coalgebra. Then there
  exists a cointegral $\gamma: C \to B$. For any  non-zero $B$-submodule
  coalgebra $D$ of $C$, $\gamma|_D: D \to B$ is clearly a
  cointegral of $D$. Hence $D$ is also $B$-cocleft.

  (R2a) Let $C$ be a $B$-module coalgebra and $\{C_i\}$ a
  chain of $B$-cocleft subcoalgebras of $C$. Since $C_i$ is an $\EE$-coalgebra for all $i$,
  by Theorem \ref{t:eclass} and (R2), $D=\bigcup_i C_i=\sum_i C_i$ is also an $\EE$-coalgebra.
  Note that ${\ol D} =
  \bigcup_i \eta_{D}(C_i)$. Since every simple subcoalgebra $E$ of
  $\ol D$ is finite-dimensional, $E$ is contained in $\eta_D(C_i)$ for some $i$.
  By Lemma \ref{l:4.2}, $D$ is
  $B$-cocleft.

  (R3) Let $C_1, C_2$ be $B$-cocleft subcoalgebras of a right $B$-module coalgebra $C$ and $D'=C_1 \land C_2$.
       By Theorem \ref{t:eclass}, $D'$ is an $\EE$-coalgebra. By \cite[Lemma 9.1.3]{Sw69},
       $$
       \ol {D'} = \eta_{D'} (D')\subseteq \eta_{D'}(C_1) \land \eta_{D'}(C_2) \subseteq \ol {D'}\,.
       $$
       For any simple coalgebra $E$ of $\ol{D'}$, $E \subseteq \eta_{D'}(C_1)$ or
       $E \subseteq \eta_{D'}(C_2)$. Hence, by Lemma \ref{l:4.2}, $D'$ is $B$-cocleft.
\end{proof}
\begin{cor}\label{c:cocleft_pro}
Let $B$ be a bialgebra.  Then every  right $B$-module coalgebra $C$
admits a unique maximal $B$-cocleft subcoalgebra $\CC(C)$ of $C$, and $\CC(C) = \CC(C) \land \CC(C)$.
Moreover, the following statements about a  right $B$-module coalgebra $C$ are equivalent :
\begin{enumerate}
\item[\rm (a)] $C$ is $B$-cocleft;
\item[\rm (b)] there exists a $B$-cocleft subcoalgebra of $C$ which contains the coradical of $C$;
\item[\rm (c)] $EB$ is $B$-cocleft for every simple subcoalgebra $E$ of $C$.
\end{enumerate}
\end{cor}
\begin{proof}
By Theorem \ref{t:cocleft}, the collection of all $B$-cocleft coalgebras is an idempotent radical class of
$\Coa{B}$. The corollary follows immediately from Proposition \ref{c-radical2}.
\end{proof}
In closing of this paper, we give another proof for a normal basis
theorem of Schneider \cite[Theorem III]{Schn90} without assuming
bijective antipode.
\begin{cor}
  Let $H$ be a Hopf algebra and $C$ a right $H$-module coalgebra. Suppose that $can : C \otimes H \to C
  \square_{\ol C} C$ is injective and $G(C)H$ contains the coradical of $C$, where $G(C)$ denotes the
  set of all group-like elements of $C$. Then $C$ is  $H$-cocleft.
\end{cor}
\begin{proof}
For $g \in G(C)$, the injectivity of $can$ implies that $\phi: H \to gH$, $h \mapsto gh$ is an $H$-module
isomorphism. Let $\psi: g H \to H$ be the inverse of $\phi$, and $\ol \psi = S \circ \psi$, where $S$ is the
antipode of $H$. Then
$$
\psi * \ol\psi (gh) = \sum h_1 S(h_2) = \e(h) 1_H = \sum S(h_1) h_2 = \ol\psi * \psi (gh)
$$
for $h \in H$. Therefore, $\psi$ is a cointegral of $gH$. In
particular, $gH$ is $H$-cocleft. By Theorem \ref{t:cocleft} and (R2), $G(C)H
= \sum_{g \in G(C)} gH$ is  $H$-cocleft. Since $G(C)H$ contains
the coradical of $C$, it follows from  Corollary \ref{c:cocleft_pro} that $C$ is
$H$-cocleft.
\end{proof}
\bibliographystyle{amsplain}

\begin{thebibliography}{1}

\bibitem{Doi83}
Yukio Doi, \emph{On the structure of relative {H}opf modules}, Comm. Algebra
  \textbf{11} (1983), no.~3, 243--255.

\bibitem{Kapl75}
Irving Kaplansky, \emph{Bialgebras}, Department of Mathematics, University of
  Chicago, Chicago, Ill., 1975, Lecture Notes in Mathematics.

\bibitem{MaDo92}
Akira Masuoka and Yukio Doi, \emph{Generalization of cleft comodule algebras},
  Comm. Algebra \textbf{20} (1992), no.~12, 3703--3721.

\bibitem{Mont93bk}
Susan Montgomery, \emph{Hopf algebras and their actions on rings}, CBMS
  Regional Conference Series in Mathematics, vol.~82, Published for the
  Conference Board of the Mathematical Sciences, Washington, DC, 1993.

\bibitem{Ng98}
Siu-Hung Ng, \emph{On the projectivity of module coalgebras}, Proc. Amer. Math.
  Soc. \textbf{126} (1998), no.~11, 3191--3198. \MR{1469428 (99a:16035)}

\bibitem{Radf77}
David~E. Radford, \emph{Pointed {H}opf algebras are free over {H}opf
  subalgebras}, J. Algebra \textbf{45} (1977), no.~2, 266--273.

\bibitem{Schn90}
Hans-J{\"u}rgen Schneider, \emph{Principal homogeneous spaces for arbitrary
  {H}opf algebras}, Israel J. Math. \textbf{72} (1990), no.~1-2, 167--195, Hopf
  algebras.

\bibitem{Sw69}
Moss~E. Sweedler, \emph{Hopf algebras}, W. A. Benjamin, Inc., New York, 1969,
  Mathematics Lecture Note Series.

\bibitem{Weibhobk}
Charles~A. Weibel, \emph{{An introduction to homological algebra}}, Cambridge
  Studies in Advanced Mathematics, no.~38, Cambridge University Press,
  Cambridge, 1994.

\end{thebibliography}
\providecommand{\bysame}{\leavevmode\hbox to3em{\hrulefill}\thinspace}
\providecommand{\MR}{\relax\ifhmode\unskip\space\fi MR }
\providecommand{\MRhref}[2]{%
  \href{http://www.ams.org/mathscinet-getitem?mr=#1}{#2}
}
\providecommand{\href}[2]{#2}

\end{document}